\RequirePackage{fix-cm}
\pdfoutput=1
\documentclass[leqno, fleqn, centertags, 12pt]{article}

\usepackage{latexsym}
\usepackage{amsmath}

\usepackage[T1]{fontenc}
\usepackage[upright]{fourier}

\usepackage{fullpage}

\newcommand{\dis}{\displaystyle}

\def\ttff{{\hspace*{0.15em}--\hspace*{0.15em}}}

\newcommand{\toto}{\longrightarrow}

\def\sss{\hspace{0.05em}\ }
\def\dss{\hspace{0.1em}\ }
\def\trs{\hspace{0.15em}\ }
\def\qss{\hspace{0.2em}\ }
\def\pss{\hspace{0.3em}\ }
\def\oss{\hspace{0.4em}\ }

\def\halfff{\hspace*{0.025em}}
\def\fff{\hspace*{0.05em}}
\def\dff{\hspace*{0.1em}}
\def\trf{\hspace*{0.15em}}
\def\qff{\hspace*{0.2em}}
\def\pff{\hspace*{0.3em}}
\def\off{\hspace*{0.4em}}

\newcommand{\nnsp}{\hspace*{-0.15em}}

\newcommand{\dnsp}{\hspace*{-0.2em}}

\makeatletter
\renewcommand{\@makefntext}[1]{
\parindent=0em
\hspace*{-0.4em}
\hbox to 0.4em{\hss\@makefnmark}\hspace*{0.4em}{#1}
}
\makeatother

\newcommand{\proof}{\vspace{0.5\bigskipamount}{\textbf{{\emph{Proof}.}}\hspace*{0.7em}}}
\newcommand{\prooftitle}[1]{\vspace{\medskipamount}{\textbf{{\emph{#1}.}}\hspace{0.7em}}}
\newcommand{\eproof}{ $\blacksquare$}
\newcommand{\esubproof}{ $\square$}

\newcounter{mysectionnumber}
\newcounter{myparnum}[mysectionnumber]
\newcounter{mylemmanum}[myparnum]
\newcommand{\myuppar}[1]{\vspace{\medskipamount}\textbf{#1}\hspace*{0.5em}}

\numberwithin{equation}{section}

\title{The\qss de\qss Bruijn--Erd\"{o}s\qss theorem\qss
in\qss incidence\qss geometry\qss
via\qss Ph.\qss Hall's\qss marriage\qss theorem}
\author{Nikolai\dss V.\dss Ivanov}
\date{}

\begin{document}

\setlength{\baselineskip}{12pt plus 0pt minus 0pt}
\setlength{\parskip}{12pt plus 0pt minus 0pt}
\setlength{\abovedisplayskip}{12pt plus 0pt minus 0pt}
\setlength{\belowdisplayskip}{12pt plus 0pt minus 0pt}

\newskip\smallskipamount \smallskipamount=3pt plus 0pt minus 0pt
\newskip\medskipamount   \medskipamount  =6pt plus 0pt minus 0pt
\newskip\bigskipamount   \bigskipamount =12pt plus 0pt minus 0pt

\maketitle

\footnotetext{\hspace*{-0.65em}\copyright\ 
Nikolai\dss V.\qss Ivanov,\oss 2017.\oss 
Neither the work reported in this paper,\qss 
nor its preparation were 
supported by any governmental 
or non-governmental agency,\qss 
foundation,\qss 
or institution.\oss}
\footnotetext{\hspace*{-0.65em}The author is grateful to F.\qss Petrov and\dss M.\qss Prokhorova for stimulating correspondence.}

\renewcommand{\baselinestretch}{1}
\selectfont

\vspace*{12pt}
\myuppar{Theorem\qss (de Bruijn--Erd\"{o}s).}
\emph{Let\qss $P$\dss 
be a finite set consisting of $n$ elements,\oss 
and let $\mathcal{L}$ be a collection of 
$m$ proper subsets of\qss $P$\nnsp.\oss  
Suppose that each\qss $l\qff \in\qff \mathcal{L}$\dss contains
at least two elements and
for every\qss $x\fff,\pff y\qff \in\qff P$\nnsp,\qss
$x\qff \neq\qff y$\qss
there exists exactly one subset\qss $l\qff \in\qff \mathcal{L}$\dss such that 
$x\fff,\pff y\qff \in\qff l$\nnsp.\oss 
Then\qss $m\qff \geqslant\qff n$\nnsp,\oss 
and\qss $m\qff =\qff n$\qss if and only if one of the following two cases
occurs:}\oss

({\fff}i\fff)\phantom{i}\qss \emph{there exists\dss $z\qff \in\qff P$\dss
such that $\mathcal{L}$ consists of\pss $P\dff \smallsetminus \{\dff z \dff\}$\qss
and all the pairs\dss $\{\dff x\fff,\pff z \dff\}$\dss with\dss $x\qff \neq\qff z$\nnsp;\oss}

\textup{({\fff}ii\fff)}\qss \emph{there exists a natural number $k$ such that
$n\qff =\qff k\dff(k\qff -\qff 1)\qff +\qff 1$\nnsp,\pss
every element of\qss $\mathcal{L}$\qss is a set consisting of $k$ elements,\oss 
every element of\pss $P$\dss belongs to exactly $k$ elements of\qss $\mathcal{L}$\nnsp,\oss
and the intersection of every two elements of\qss $\mathcal{L}$\qss is non-empty.}

\vspace*{6pt}  
The elements of $P$ are called \emph{points}\pss
and the elements of $\mathcal{L}$ \emph{lines}.\qss
For a point $a$ let $k_{\fff a}$ be the number of lines containing $a$\nnsp,\qss
and for a line $l$\dss let\dss $s_{\fff l}$\dss be the number of points in $l$\nnsp,\qss
i.e. the cardinality of the set $l$\nnsp.\qss
Following de Bruijn--Erd\"{o}s\dss \cite{db-e},\qss we start with two observations.\qss
First\halfff,\qss counting the pairs 
$(a\fff,\pff l\dff)\qff \in\qff
P\times\mathcal{L}$\dss such that\dss
$a\qff \in\qff l$\qss in two ways leads to the equality
\begin{equation}
\label{sums}
\quad
\sum_{\qff a\dff \in\dff P}\qff k_{\fff a}
\off\off =\off\off
\sum_{l\dff \in\dff \mathcal{L}}\qff s_{\fff l}\qff. 
\end{equation}
Second,\qss if
$(a\fff,\pff l\dff)\qff \in\qff
P\times\mathcal{L}$\dss and\dss
$a\qff \not\in\qff l$\nnsp,\qss then for every $x\qff \in\qff l$\dss there
is a line containing $\{\fff a\fff,\pff x \dff\}$
and these lines are distinct\halfff.\oss
It follows that in this case
\begin{equation}
\label{inequality}
\quad
s_{\fff l}
\off \leqslant\off
k_{\fff a}\qff.
\end{equation}

\vspace*{-12pt}
\myuppar{Systems of distinct representatives.}
A\qss (finite)\qss family of subsets of a set $S$ is a map\qss 
$\varphi \colon I\toto 2^S$ 
from a finite set $I$
to the power set $2^S$\nnsp.\qss
Usually $\varphi\fff(i)$ is denoted by $\varphi_i$\dnsp.\qss
A\qss \emph{system of distinct representatives}\qss for $\varphi$
is an\qss \emph{injective}\qss map\qss
$f \colon I\toto S$\qss
such that\qss
$f(i)\qff \in\qff \varphi_i$\qss for all\qss 
$i\qff \in\qff I$\nnsp.\qss 
We will denote by $|\dff X\dff|$ the number of elements of a set $X$\nnsp.\qss

\myuppar{Theorem\qss (Ph.\dss Hall).}
\emph{A system of distinct representatives for $\varphi$ exists
if and only if\dss
the union
\begin{equation}
\label{union}
\quad
\bigcup\nolimits_{\dff i\dff \in\dff J}\qff \varphi_i
\end{equation}
contains\qss $\geqslant\qff |\dff J\dff|$\qss elements
for every subset\qss
$J\qff \subset\qff I$\nnsp.\oss}

\myuppar{Corollary.}
\emph{Suppose that for every non-empty subset\qss
$J\qff \subset\qff I$\qss
not equal to $I$ itself
the union\qss (\ref{union})\qss contains\qss $>\qff |\dff J\dff|$\qss elements.\qss
Then for every\qss $j\qff \in\qff I$\dss 
and every element $a\qff \in\qff \varphi_j$
there exists a system of distinct representatives $f \colon I\toto S$ for $\varphi$
such that\qss $f(j)\off =\off a$\nnsp.}

\proof\qss
Let\qss $K\qff =\qff I\dff \smallsetminus\dff \{\dff j\dff\}$\nnsp,\qss 
and let\qss $\psi_i\off =\off \varphi_i\dff \smallsetminus\dff  \{\dff a\dff\}$\qss
for every\qss $i\qff \in\qff K$\nnsp.\qss 
Then $\psi$ is a family of subsets of $S$\nnsp.\qss
The theorem applies to $\psi$ and implies that there exists
a system of distinct representatives\qss $g \colon K\toto S$\qss 
for $\psi$\nnsp.\oss
The map\qss $f \colon I\toto S$\qss defined by\qss
$f(i)\off =\off g(i)$\qss for\qss $i\qff \in\qff K$\qss
and\qss $f(j)\off =\off a$\qss is the desired system of distinct representatives 
for $\varphi$\nnsp.\oss  \eproof

\vspace*{6pt}
In fact\halfff,\qss 
a standard proof\dss of\dss the\dss Ph.\qss Hall's\sss theorem using an induction
by the number of elements of $I$ includes a proof of this Corollary.\qss
See,\qss for example,\qss the textbook\trs \cite{c}\trs by\dss P.\qss Cameron.

\myuppar{The family of complements to lines.}
Let\qss $I\qff =\qff \mathcal{L}$\qss and\qss
$\varphi_l\qff =\qff P\dff \smallsetminus\dff l$\qss for\qss $l\qff \in\qff \mathcal{L}$\dnsp.\oss
Then $\varphi$ is family of subsets of $P$\dnsp.\oss
If\qss $l\qff \in\qff \mathcal{L}$\dnsp,\oss then\qss $l\qff \neq\qff P$\qss
and hence $\varphi_l\qff \neq\qff \emptyset$\dnsp.\oss
If\qss $l\fff,\pff l'\qff \in\qff \mathcal{L}$\qss
and\qss $l\qff \neq\qff l'$\nnsp,\oss
then\qss $l\dff \cap\dff l'$ contains no more than one point and hence\qss 
$\varphi_l\dff \cup\dff \varphi_{l'}$\qss
is equal either to $P$\dnsp,\oss or to $P$ with one point removed.\oss
It follows that the same is true for the union\qss (\ref{union})\qss
if\dss $J$\dss contains at least two lines,\oss
and hence in this case\qss (\ref{union})\qss contains\qss 
$\geqslant\qff n\qff -\qff 1$\qss points.

Let $l$ be a line and let\qss
$x\fff,\pff y\qff \in\qff l$\qss and\qss $x\qff \neq\qff y$\nnsp.\qss
Let $a$ be a point not in $l$\nnsp.\qss
Then the two lines $l_x\fff,\pff l_y$ containing $a$ and $x\fff,\pff y$ respectively
intersect only in $a$\nnsp.\qss
It follows that the intersection of three lines $l\fff,\pff l_x\fff,\pff l_y$
is empty,\qss
and hence the union of their complements is equal to $P$\dnsp.\qss
It follows that for\qss $J\qff =\qff I\qff =\qff \mathcal{L}$\qss 
the union\qss (\ref{union})\qss is equal to $P$
and hence contains $n$ elements.\qss

\prooftitle{Proof of the inequality\qss $m\qff \geqslant\qff n$\nnsp}
Suppose that $m\qff \leqslant\qff n$\nnsp.\oss
Let\qss $J\qff \subset\qff \mathcal{L}$\dnsp.\oss
If\qss $|\dff J\dff|\qff =\qff 1$\nnsp,\qss
then\qss (\ref{union})\qss is the complement of a line 
and hence contains\qss $\geqslant\qff 1$\qss points.\qss 
If\qss $2\qff \leqslant\qff |\dff J\dff |\qff \leqslant\qff m\qff -\qff 1$\nnsp,\qss
then\qss (\ref{union})\qss contains\qss
$\geqslant\qff n\qff -\qff 1\qff \geqslant\qff m\qff -\qff 1$\qss points,\qss
and hence contains\qss $\geqslant \qff |\dff J\dff |$\qss elements.\qss
If\qss $|\dff J\dff|\qff =\qff m$\nnsp,\qss
then\qss (\ref{union})\qss contains\qss
$n\qff \geqslant\qff m\qff =\qff |\dff J\dff|$\qss elements.\qss
Therefore,\oss by\dss Hall's theorem there exists a system of distinct representatives
for $\varphi$\nnsp.\qss
Such a system of distinct representatives is an injective map\qss
$f \colon \mathcal{L}\toto P$\qss
such that\qss
$f(l\fff)\off \not\in\off l$\qss for every\qss
$l\qff \in\qff \mathcal{L}$\nnsp.\qss
By the inequality\qss (\ref{inequality}),\qss 
this implies that\qss
$\dis
s_l\qff \leqslant\qff k_{\fff f(l\fff)}$\qss
for every\qss
$l\qff \in\qff \mathcal{L}$\nnsp.\qss
By summing all these inequalities and using the injectivity of $f$ we see that
\begin{equation}
\label{three-sums}
\quad
\sum_{l\dff \in\dff \mathcal{L}}\qff s_{\fff l}
\off\off \leqslant\off\off
\sum_{l\dff \in\dff \mathcal{L}}\qff k_{\fff f(l\fff)}
\off\off \leqslant\off\off
\sum_{\qff a\dff \in\dff P}\qff k_{a}\qff.
\end{equation}
Moreover\halfff,\qss
the second inequality is strict unless\qss $m\qff =\qff n$\qss
(otherwise the last sum has an additional summand\qss $k_a\qff >\qff 0$\qss
compared to the previous one).\qss
But\qss (\ref{sums})\qss implies that both inequalities in\qss
(\ref{three-sums})\qss should be actually equalities.\oss
It follows that\qss $m\qff =\qff n$\nnsp.\qss  \eproof

\prooftitle{The case\qss $m\qff =\qff n$\dnsp.\oss  I\dff}
Suppose that there is a point\qss $a\qff \in\qff P$\qss
contained in\qss $m\qff -\qff 1$\qss lines.\qss
Each of these lines contains at least one point in addition to $a$\nnsp.\oss
Since\qss $m\qff =\qff n$\nnsp,\qss
there are no other points.\qss
On the other hand,\qss
there is only one more line.\qss
This line should contain all the points different from $a$\nnsp.\qss
It follows that we are in the case\qss ({\fff}i\fff)\qss of the
de Bruijn-Erd\"{o}s theorem.\qss  \esubproof

\prooftitle{The case\qss $m\qff =\qff n$\dnsp.\oss II\dff}
Suppose that\qss $P\qff =\qff l\dff \cup\dff l'$\qss
for some\qss $l\fff,\pff l'\qff \in\qff \mathcal{L}$\qss
and\qss $l\dff \cap\dff l'\qff \neq\qff \emptyset$\dnsp.\oss
Then the intersection\qss $l\dff \cap\dff l'$\qss consists of exactly one point\halfff,\oss
say,\qss the point $z$\nnsp.\qss
It follows that the number of points\qss $n\off =\off s_l\qff +\qff s_{l'}\qff -\qff 1$\nnsp.\qss
Let us count the lines.\oss

In addition to the lines\qss $l\fff,\pff l'$\dnsp,\oss
for every two points\qss 
$b\qff \in\qff l\qff \smallsetminus\dss \{\dff z\dff\}$\qss
and\qss 
$b'\qff \in\qff l'\qff \smallsetminus\dss \{\dff z\dff\}$\qss 
there is a line containing\qss $b\fff,\pff b'$\dnsp.\oss
All these lines are distinct\halfff,\qss
and hence there are\qss 
$\geqslant\qff 2\qff +\qff (s_l\qff -\qff 1)(s_{l'}\qff -\qff 1)$\qss
lines.\oss
We may assume that\qss $s_l\qff \geqslant\qff s_{l'}$\dnsp.\oss
If\qss $s_{l'}\qff \geqslant\qff 3$\dnsp,\oss
then the number $m$ of\trs lines\dss is\qss
\[
\quad
\geqslant\off
2\qff +\qff (s_l\qff -\qff 1)(s_{l'}\qff -\qff 1)
\off \geqslant\off 
2\qff +\qff 2\dff(s_l\qff -\qff 1)
\off =\off 
2\fff s_l
\off \geqslant\off 
s_l\qff +\qff s_{l'}
\off =\off 
n\qff +\qff 1\dff,
\]
contrary to the assumption\qss $m\qff =\qff n$\nnsp.\qss
It follows that\qss $s_{l'}\qff =\qff 2$\qss
and hence all points,\qss except the only point in\qss
$l'\dff \smallsetminus\dff l$\nnsp,\qss belong to the line $l$\dnsp.\oss
Let $a$ be this point\halfff,\oss so\qss
$\{\fff a\trf\}\qff =\qff l'\dff \smallsetminus\dff l$\nnsp.\oss
In this case\qss $s_l\off =\off n\qff -\qff 1\off =\off m\qff -\qff 1$\qss
and $a$ belongs to\qss $\geqslant\off s_l\off =\off m\qff -\qff 1$\qss
lines.\oss
Therefore,\qss we are in the situation of the case\qss \textbf{I}\qss again,\qss
and hence in the case\qss ({\fff}i\fff)\qss of the
de Bruijn-Erd\"{o}s theorem.\qss  \esubproof

\prooftitle{The case\qss $m\qff =\qff n$\dnsp.\oss III\dff}
Suppose that we are neither in the case\qss \textbf{I},\qss
nor in the case\qss \textbf{II}.\oss 
In particular\halfff,\qss no point is contained in\qss $m\qff -\qff 1$\qss lines.\qss 
In this case,\qss if\dss $J$\dss is a set of\qss
$m\qff -\qff 1$ lines,\qss
then\qss (\ref{union})\qss contains\qss $n\qff >\qff m\qff -\qff 1$\qss points.\qss
If\dss $J$\dss is a set of\dss $p$\dss lines and\qss
$p\qff <\qff m\qff -\qff 1$\nnsp,\qss
then\qss (\ref{union})\qss contains\qss 
$\geqslant\qff n\qff -\qff 1\qff =\qff m\qff -\qff 1\qff >\qff p$\qss points.\qss
If\qss $p\qff =\qff 1$\qss and\qss (\ref{union})\qss 
contains only $1$ point\halfff,\qss then this point is contained in\qss $m\qff -\qff 1$\qss
lines and we are in the case\qss \textbf{I},\qss
contrary to the assumption.\qss
Therefore,\qss if\qss $p\qff =\qff 1$\nnsp,\qss then\qss (\ref{union})\qss 
contains\qss $\geqslant\qff 2\qff >\qff p$\qss points.\qss
It follows that in this case the Corollary applies.

Let $f$ be a system of distinct representatives for $\varphi$\nnsp.\qss
Since $f$ is injective and\qss $m\qff =\qff n$\nnsp,\qss
the map $f$ is actually a bijection.\qss
Therefore,\qss the inequalities\qss
$\dis
s_l\off \leqslant\off k_{\fff f(l\fff)}$\qss
together with\qss (\ref{sums})\qss imply that\qss
$s_l\off =\off k_{\fff f(l\fff)}$\qss
for every\qss
$l\qff \in\qff \mathcal{L}$\nnsp.\oss
By the Corollary,\qss
given a line $l$ and a point\qss $a\qff \not\in\qff l$\nnsp,\qss
one can choose $f$ is such a way that\qss $f(l\fff)\qff =\qff a$\nnsp.\oss
It follows that\qss $s_l\off =\off k_a$\qss
if\qss $a\qff \not\in\qff l$\nnsp.\oss
This implies,\qss
in particular\halfff,\qss
that\dss if\qss $a\qff \not\in\qff l$\nnsp,\qss then
every line containing $a$ intersects $l$\dnsp.\oss
In turn,\qss this implies that every two lines intersect\halfff.\oss

Suppose that\qss $l\fff,\pff l'\qff \in\qff \mathcal{L}$\qss and\qss $l\qff \neq\qff l'$\dnsp.\oss
By the previous paragraph\qss $l\dff \cap\dff l'\qff \neq\qff \emptyset$\qss
and hence if\qss $P\qff =\qff l\dff \cup\dff l'$\dnsp,\qss
then we are in the case\qss \textbf{II},\qss
contrary to the assumption.\oss
Therefore\qss $P\qff \neq\qff l\dff \cup\dff l'$\dnsp.\oss
If\qss $a\qff \not\in\qff l\dff \cup\dff l'$\dnsp,\oss
then\qss $s_l\qff =\qff k_a\qff =\qff s_{l'}$\nnsp.\oss
It follows that all the numbers\qss 
$s_l\fff,\pff k_a$\qss are equal.\qss
Let $k$ be the common value of these numbers\qss $s_l\fff,\pff k_a$\nnsp.\oss
Let\qss $b\qff \in\qff P$\dnsp.\oss
There are $k$ lines containing $b$ and each of the\qss $n\qff -\qff 1$\qss 
points not equal to $b$ is contained in 
one and only one of these lines.\oss
Moreover\halfff,\qss each of these lines consists of $b$ and exactly\dss $k\qff -\qff 1$\dss points\qss
not equal to $b$\nnsp.\oss
It follows that\qss $n\off =\off k\dff(k\qff -\qff 1)\qff +\qff 1$\dnsp.\oss
Since by the previous paragraph every two lines intersect\halfff,\qss
we are in the case\qss ({\fff}ii\fff)\qss of the
de Bruijn-Erd\"{o}s theorem.\oss  \esubproof  \eproof

\begin{flushright}

April\qss 17,\oss 2017
 
https\halfff:/\!/\hspace*{-0.06em}nikolaivivanov.com

E-mail\halfff:\oss nikolai.v.ivanov{\fff}@{\dff}icloud.com

\end{flushright}

\end{document}